\documentclass{article}

\usepackage[utf8]{inputenc}
\usepackage{amsmath, amssymb, amsthm, bm, mathrsfs}
\usepackage[top=1.5cm, bottom=1.5cm, left=2cm, right=2cm]{geometry}
\usepackage{csquotes, authblk, cases, enumitem}
\usepackage{graphicx, float, caption, subcaption, adjustbox}
\usepackage[hidelinks]{hyperref}
\usepackage{cleveref}
\usepackage[numbers,sort]{natbib}
\usepackage[section]{algorithm}
\usepackage{algpseudocode}
\usepackage{array, tabularx, xcolor, colortbl}
\usepackage{color, xcolor}
\usepackage[para]{footmisc}
\usepackage[normalem]{ulem}
\counterwithin{figure}{section}
\counterwithin{table}{section}
\numberwithin{equation}{section}
\graphicspath{{Figures/}}
\definecolor{darkred}{rgb}{0.8,0,0}

\newenvironment{frcseries}{\fontfamily{frc}\selectfont}{}

\theoremstyle{definition}

\title{ Direct Element-Level Mass Matrix Inversion }
\author{Eli Hanukah \thanks{eli.hanukah@braude.ac.il}}
\affil{\small Braude College of Engineering, Karmiel. P.O. Box 78, 2161002, Israel}
\date{\today}

\begin{document}

\maketitle

\begin{abstract}
Explicit finite element methods for transient solid dynamics commonly employ lumped mass matrices to avoid the computational expense of solving systems of equations for nodal accelerations at each time step. However, diagonal mass matrices can inadequately represent inertia, particularly for higher-order elements and in coarse or distorted-mesh regions. Consequently, many implementations have limited the use of higher-order elements in explicit dynamic analyses, relying predominantly on linear elements and lumped mass formulations. This paper presents a novel, systematic methodology for directly evaluating the inverse of non-diagonal element-level mass matrices for solid isoparametric elements of arbitrary interpolation order. The method applies to a family of mixed mass matrices defined as linear combinations of the lumped and consistent mass matrices, with a weighting parameter governing the transition between the diagonal, consistent, and intermediate formulations. By enabling the direct evaluation of element-level mass matrix inverses, the proposed methodology provides a general framework for incorporating non-diagonal mass formulations into explicit finite element algorithms, including those employing higher-order elements, such as quadratic elements. To illustrate the approach, we derive an explicit formula for the element-level inverse of the mixed mass matrix for the eight-node hexahedral element. The resulting expression provides a computationally efficient way to evaluate the element-level inverse without solving the global system. The proposed methodology can improve inertial representation in transient dynamic analyses, particularly with higher-order elements and in coarse or distorted-mesh regions, while retaining the computational advantages of explicit time integration. It is amenable to implementation in element-by-element computational frameworks commonly used in explicit finite element analysis.\\

\end{abstract}

\section{Introduction}
\label{Sec. Introduction}

\noindent \indent The finite element method is widely used to numerically analyze transient dynamic problems in solid mechanics, including impact, wave propagation, and structural response under rapidly varying loads. In these applications, the spatial discretization of the equations of motion leads to a system of second-order differential equations of the form 
\begin{equation}
    \mathbf{M}\ddot{\mathbf{u}}\left(t\right) = \mathbf{f}^{ext}(\mathbf{u}(t)) - \mathbf{f}^{int}(\mathbf{u}(t))
\nonumber
\end{equation}
where \(\mathbf{M} \) denotes the mass matrix, \( t \) is time, \( \mathbf{u} \) is the vector of nodal displacements, \( \mathbf{f}^{int} \) and \( \mathbf{f}^{ext} \) are the internal and external force vectors, respectively. Representing inertia through the mass matrix is an important part of spatial discretization and can affect the accuracy and computational cost of the resulting dynamic solution. \\
\noindent \indent Two principal approaches to mass matrix construction are commonly employed: the consistent mass formulation, which arises naturally from the isoparametric approximation, and the lumped mass formulation, which redistributes mass to produce a diagonal matrix. The consistent mass matrix retains coupling between nodal degrees of freedom associated with the finite element interpolation, whereas the diagonal lumped mass matrix, due to its trivial inverse, allows independent evaluation of nodal accelerations. Researchers have extensively investigated the properties and accuracy of these formulations in structural dynamics and transient elasticity ~\cite{Tong1971, Wu2006}. \\
\noindent \indent Explicit time integration methods are particularly attractive for nonlinear transient analyses involving large deformations, contact, impact, and complex constitutive behavior. Their computational efficiency depends strongly on how the mass matrix is treated. In particular, a diagonal lumped mass matrix allows the acceleration vector to be computed through independent nodal operations, avoiding solving a coupled global system of equations at each time step. Consequently, lumped mass formulations are widely employed in explicit finite element analysis. Nevertheless, the advantages of mass lumping must be weighed against its effects on the discrete inertia representation and the accuracy of the computed dynamic response ~\cite{Wu2006, Yang2017}. \\
\noindent \indent Although diagonal mass matrices offer substantial computational advantages, retaining off-diagonal mass terms may be desirable when a more faithful representation of the discretized system's inertial properties is required. This consideration is particularly relevant to coarse-mesh and distorted-mesh simulations, where the spatial discretization's influence on the dynamic response may be more pronounced. \\
\noindent \indent Mixed mass formulations, constructed as linear combinations of the lumped and consistent mass matrices, provide a systematic way to interpolate between these two representations. A weighting parameter controls each formulation's contribution, allowing the diagonal, consistent, and intermediate mass matrices to be treated within a unified framework. However, using such non-diagonal formulations in explicit dynamics requires an efficient way to handle the associated mass matrix in the acceleration update ~\cite{lim2008nonlinear}. \\
\noindent \indent Standard implementations use diagonal lumped mass matrices; this approach often degrades spatial accuracy and fails for higher-order finite elements because it can produce non-physical zero or negative nodal masses. To address these limitations, recent paradigms have shifted toward directly constructing consistent, sparse inverse mass matrices using multi-field variational principles. For instance, Tkachuk and Bischoff ~\cite{Tkachuk2013} developed a variational formulation based on a modified Hamilton’s principle that incorporates biorthogonal momentum fields. This strategy directly generates a sparse reciprocal mass matrix at the element level, eliminating the need for global inversion and proving highly effective for selective mass scaling in low-order, coarse, or highly constrained meshes. To address the challenges of higher-order finite element formulations, Kolman et al. \cite{Kolman2020} extended this methodology by introducing a penalized three-field formulation using density-weighted dual shape functions. Their approach accurately preserves total system mass and captures complex high-frequency linear free-vibration spectra, bypassing the traditional mathematical pitfalls of mass lumping in higher-order elements by directly assembling a sparse global inverse matrix through localized, element-level dual-basis projections. Moving beyond finite element frameworks, Nairn ~\cite{Nairn2026} expanded the use of the non-diagonal full mass matrix within particle-based methods. This approach utilizes \(\text{FMPM}(k)\) order-based Taylor approximations—modeled via Neumann series expansions—to efficiently approximate the global inverse through sparse matrix-vector multiplications, while adapting the formulation to handle explicit boundary conditions and contact mechanics in Material Point Method (MPM) simulations. \\
\noindent \indent To preserve the computational efficiency of explicit time-stepping schemes without sacrificing the accuracy of consistent formulations, researchers have widely explored element-by-element (EBE) localized inversion and direct integration strategies. In their foundational work, ~\cite{hughes1978functional} introduced a functional EBE framework that bypasses global structural matrix assembly by decomposing the transient equations into localized element operators. By performing operations and inversions directly at the element level, the approach retains non-diagonal mass couplings while maintaining a highly parallelizable explicit update structure. Building on these localized strategies, ~\cite{karakostas2017improving} optimized mass and inverse-mass matrix computation specifically for explicit dynamics. Their approach streamlines element-level inverse calculations while reducing the required Gauss integration points, effectively mitigating the dense global fill-in penalty traditionally associated with non-diagonal mass formulations. While these element-level techniques capture high-frequency wave mechanics on coarse meshes by retaining off-diagonal mass coupling, they remain highly sensitive to geometric element distortion, which can degrade local Jacobian-matrix conditioning and compromise the explicit solver's stability. Building on these efficiency goals for specific spatial discretizations, ~\cite{hanukah2018improving} introduced a semi-analytical approach for 8-node and 20-node hexahedral elements. Their scheme drastically reduces computational overhead by using optimized, low-point integration rules (such as a 4-point rule replacing standard 8-point or 14-point quadratures) to evaluate element mass matrices, while also introducing a novel approach to directly compute the inverse mass matrix of 8-node hexahedral elements. \\
\noindent \indent This study develops a systematic method for directly evaluating the inverse of element-level mixed mass matrices for standard solid isoparametric elements in two and three dimensions, including both lower- and higher-order formulations. Each entry of the general inverse is expressed as a polynomial in Jacobian-related variables divided by another polynomial in the same variables. The coefficients are determined analytically. The paper establishes the convergence of the proposed approximation to the exact inverse and provides a systematic approach for improving its accuracy by successively adding terms. \\
\noindent \indent In addition, we derive a detailed expression for the inverse of the mixed mass matrix associated with the widely used eight-node hexahedral element and establish its mathematical properties. In particular, the proposed expression is positive definite, provided that the Jacobian determinant evaluated at the element centroid is positive. In terms of computational efficiency, the proposed expression requires fewer operations than conventional one-point numerical integration for mass matrix evaluation and eliminates the need for the matrix inversion algorithm associated with the latter approach. \\
\noindent \indent The proposed method is intended to facilitate incorporating non-diagonal mass formulations into explicit transient finite element algorithms through element-by-element computation. This approach may be particularly beneficial in coarse- and distorted-mesh simulations, where retaining off-diagonal mass terms may improve inertia representation. The method also allows selective application to specific elements or uniform application across the mesh, offering flexibility in balancing computational cost and solution accuracy. The practical effects of this approach on complete transient dynamic simulations remain an important subject for numerical assessment. \\
\noindent \indent The remainder of the paper is organized as follows. First, in Section \ref{Sec. background}, we recall mathematical definitions associated with the consistent, lumped (diagonal), and mixed mass matrices. Then, we recall the Standard (ST), the Semi-Analytical (SA), and the J-Derivatives (JD) approaches to mass matrix integration and representation. Second, in Section \ref{Sec. DMI}, we derive the general mathematical formulation for the direct inverse of the element-level mixed mass matrix. Next, in Section \ref{Sec. RGM}, we develop a practical, ready-to-implement formula for the widely-used eight-node hexahedral element. We demonstrate its efficiency and positive definiteness, even for highly distorted meshes. Finally, in Section \ref{Sec. conclusions}, we present the principal conclusions and directions for future research. \\

\section{Mathematical background}
\label{Sec. background}

\noindent \indent We begin by introducing the notation and recalling the definitions required for the subsequent formulation. The number of nodes per element is denoted by \(n_{\mathrm{nodes}}\), while the number of spatial dimensions, corresponding to either a two- or three-dimensional formulation, is denoted by \(n_{\mathrm{dim}}\). Throughout this study, uppercase indices range over the element nodes, \( (I, J, K=1,\ldots,n_{\mathrm{nodes}}) \), whereas lowercase indices range over the spatial dimensions, \( (m,n,r=1,\ldots,n_{\mathrm{dim}}) \). \\
\noindent \indent For a purely mechanical formulation, the number of element degrees of freedom is
\(n_{\mathrm{dof}}=n_{\mathrm{nodes}} \times n_{\mathrm{dim}}\). In the nodal representation, the element mass matrix is denoted by \(M_{IJ}\) and has dimensions \(n_{\mathrm{nodes}}\times n_{\mathrm{nodes}}\). In the scalar degree-of-freedom representation, the corresponding matrix is denoted by \(M_{ij}\) and has dimensions \(n_{\mathrm{dof}}\times n_{\mathrm{dof}}\). For the standard mechanical element-level mass matrix, these representations are related by a block-diagonal expansion in which each nodal entry is repeated along the corresponding diagonal component blocks, while the off-diagonal component blocks vanish. Thus, the two representations contain the same independent scalar entries. The same convention applies to the inverse of the element-level mass matrix, denoted by \(W_{IJ}\) or \(W_{ij}\), respectively. We use the letter \(W\) because it visually resembles the letter \(M\) inverted. \\
\noindent \indent The transformation from the nodal representation to the scalar degree-of-freedom representation is defined by the index rule in Eq.~(\ref{eq: transform}). In what follows, we adopt the compact nodal notation \(M_{IJ}\) and \(W_{IJ}\), without loss of generality for the block-diagonal mechanical mass matrices considered here.
\begin{equation}
    \label{eq: transform}
    \begin{array}{c} 
    A_{i,j} 
    \\  
    \left( i,j=1,\ldots,n_{dof} \right) 
    \end{array}=
           \left\{\begin{array}{c}
                 \\
                 A_{IJ}
                 \\ \\ \\ \\
                 0  
                 \end{array}  \qquad  \begin{array}{c}
                 i=n_{dim}(I-1)+m \\
                 j=n_{dim}(J-1)+m \\
                 \left( I,J=1,\ldots,n_{nodes},m=1,\ldots,n_{dim} \right) \\ \\ \\
                 otherwise  
                 \end{array}\right\}
\end{equation}
\noindent \indent The shape functions are denoted by \(N_I\), and the initial mass density is denoted by \(\rho_0\). The differential volume in the initial configuration is related to the differential volume in the parent element by \( dV=J\,d\Box \), where \(J=\det(\bm{J})\) is the determinant of the Jacobian matrix and \(d\Box\) denotes the differential volume in the parent domain. Specifically,
\[
d\Box =
\begin{cases}
d\xi\,d\eta, & n_{\mathrm{dim}}=2,\\
d\xi\,d\eta\,d\zeta, & n_{\mathrm{dim}}=3,
\end{cases}
\]
where \(\bm{\xi}=(\xi,\eta)\) or \(\bm{\xi}=(\xi,\eta,\zeta)\) denotes the local coordinates, as appropriate. Accordingly, integration of an arbitrary scalar function \(A\) over the element domain \(V_e\) can be expressed as
\[ \int_{V_e} A\,dV=\int_{V_{\Box}} A\,J\,d\Box \]
For example, for a three-dimensional hexahedral element whose parent domain is the unit cube, the integral becomes
\[ \int_{V_{\Box}} A\,J\,d\Box = \int_{-1}^{1}\int_{-1}^{1}\int_{-1}^{1} A\,J\,d\xi\,d\eta\,d\zeta\]
\noindent \indent Following the isoparametric formulation, the position vector of an arbitrary material point in the initial configuration is interpolated as
\[ \bm{X}(\bm{\xi}) = 
\sum_{I=1}^{n_{\mathrm{nodes}}}
N_I(\bm{\xi})\,\bm{X}_I \]
where \(\bm{X}_I\) denotes the initial position of node \(I\). The nodal positions are expressed in the initial global basis
\((\bm{E}_1,\ldots,\bm{E}_{n_{\mathrm{dim}}})\) as
\[
\bm{X}_I
=
\sum_{m=1}^{n_{\mathrm{dim}}}
X_{Im}\bm{E}_m
\]
Consequently, the global Cartesian components of the position vector are
\[
X_m=\bm{X}\cdot\bm{E}_m
=\sum_{I=1}^{n_{\mathrm{nodes}}}N_I X_{Im}
\]

The Jacobian matrix is defined as the gradient of the position vector with respect to the local coordinates ~\( \bm{J}=Grad_{\bm{\xi}} \left( \bm{X} \right) \) (see, e.g., \cite{wriggers2010nonlinear}). Its components are
\[
J_{mn}=X_{m;n}
=\sum_{I=1}^{n_{\mathrm{nodes}}}
N_{I;n}(\bm{\xi})X_{Im},
\qquad
(m,n=1,\ldots,n_{\mathrm{dim}})
\]
where a semicolon denotes partial differentiation with respect to a local coordinate. For example, \( N_{I;2}=\frac{\partial N_I}{\partial\eta} \).

The Jacobian determinant depends on both the local coordinates \( \bm{\xi} \) and the element geometry \( X_{Km} \). The shape functions \(N_I(\bm{\xi})\) are polynomials in the local coordinates and do not depend directly on the nodal positions. By contrast, the nodal coordinates \(X_{Kn}\) are fixed parameters for a given element geometry. If the shape functions have polynomial degree \(\kappa\), the components of the Jacobian matrix have degree at most \( (\kappa-1) \) in the local coordinates. Hence, the Jacobian determinant \( J(\bm{\xi};X_{Kn})=\det(\bm{J}) \) is a geometry-dependent polynomial of degree at most \(n_{\mathrm{dim}}(\kappa-1)\). For example, in a ten-node tetrahedral element, the shape functions are complete order quadratic polynomials. The Jacobian matrix components are therefore linear in the local coordinates, and its determinant is a geometry-dependent polynomial of third degree.
\noindent \indent The Jacobian determinant can be represented or approximated in several ways. One approach is to expand it in a multivariate Taylor series about the element centroid, \( \bm{\xi}^*=\left( \xi^*,\eta^*,\zeta^* \right) \). This representation forms the basis of the J-Derivatives (JD) mass-matrix integration method \cite{hanukah_2025_element,hanukah_2025_JD10tet}.
\begin{equation}
    \begin{gathered}
    J^{JD}=J{\Big|}_{\bm{\xi}^{*}}+
           \left( \xi-\xi^{*} \right) \frac{\partial{J}}{\partial{\xi}} {\Big|}_{\bm{\xi}^{*}}+
           \left( \eta-\eta^{*} \right) \frac{\partial{J}}{\partial{\eta}} {\Big|}_{\bm{\xi}^{*}}+
           \left( \zeta-\zeta^{*} \right) \frac{\partial{J}}{\partial{\zeta}} {\Big|}_{\bm{\xi}^{*}}+...
    \nonumber
    \end{gathered}
\end{equation}
The expansion is truncated after a finite number of terms, yielding the approximation \( J^{JD} \approx J \). The number of terms \( n \) depends on the element type and the polynomial structure of the Jacobian determinant. The resulting representation separates dependence on local coordinates from dependence on the geometry (nodal positions/coordinates), as expressed in Eq.~(\ref{eq: J_JD}).
\begin{equation}
    \label{eq: J_JD}
    \begin{gathered}
    J^{JD}= 
    \sum_{q=0}^{n-1} \tilde{N}_q \left( \bm{\xi} \right) \tilde{J}_q \left( X_{Kn} \right), \qquad  
    \left( K=1,\ldots,n_{nodes}, n=1,\ldots,n_{dim} \right)
    \end{gathered}
\end{equation}
\begin{equation}
    \label{eq: N_tilde}
    \begin{gathered}
    \tilde{N_0}=1, \qquad
    \tilde{N_1}= \left( \xi-\xi^* \right),\qquad
    \tilde{N_2}= \left( \eta-\eta^* \right),\qquad
    \tilde{N_3}= \left( \zeta-\zeta^* \right),\qquad ...
    \end{gathered}
\end{equation}
\begin{equation}
    \label{eq: J_tilde}
    \begin{gathered}
    \tilde{J}_0 \left( X_{Kn} \right) = J {\Big|}_{\bm{\xi}^*}, \qquad
    \tilde{J_1} \left( X_{Kn} \right) = 
    \frac{\partial{J}}{\partial{\xi}} {\Big|}_{\bm{\xi}^*}, \qquad
    \tilde{J}_2 \left( X_{Kn} \right) = \frac{\partial{J}}{\partial{\eta}} {\Big|}_{\bm{\xi}^*}, \qquad
    \tilde{J}_3 \left( X_{Kn} \right) = \frac{\partial{J}}{\partial{\zeta}} {\Big|}_{\bm{\xi}^*}, \qquad ...
    \end{gathered}
\end{equation}
\noindent \indent The Semi-Analytical (SA) approach represents or approximates the Jacobian determinant using polynomial interpolation functions  \( \hat{N}_q \left( \bm{\xi} \right) \), analogous to the isoparametric interpolation of the geometry. These functions satisfy the interpolation property ~\( \hat{N}_q {\Big|}_{point \; w} = \delta_{qw} \), (see, e.g., \cite{hanukah2014semi,hanukah2014consistent,hanukah2014exact,hanukah2014mass,hanukah2015semi, HG2_2016_reduce_int_point}). The resulting representation is given by Eq.~(\ref{eq: J_SA}), where \( J_q={J} {\Big|}_{point \; q} \) denotes the Jacobian determinant evaluated at sampling point. As in the JD formulation, this representation separates the coordinate-dependent interpolation functions from the geometry-dependent coefficients.
\begin{equation}
    \label{eq: J_SA}
    \begin{gathered}
    J^{SA}= 
    \sum_{q=1}^{n} \hat{N}_q \left( \bm{\xi} \right) J_q \left( X_{Kn} \right), \qquad  
    \left( K=1,\ldots,n_{nodes}, n=1,\ldots,n_{dim} \right)
    \end{gathered}
\end{equation}
Separation of variables between the coordinate-dependent polynomials \( \hat{N}_q \left( \bm{\xi} \right) \) and mesh-dependent \( J_q \left( X_{Kn} \right) \) is maintained. \\
\noindent \indent The consistent element mass matrix, denoted by \( M_{IJ}^{cons} \), \( \left( I, J=1,\ldots,n_{nodes} \right) \), is defined as the integral of the product of the shape functions weighted by the initial mass density over the element domain.
\begin{equation}
    \label{eq: cmm_def}
    \begin{gathered}
    M_{IJ}^{cons}=\int_{V_e} \rho_0 N_I N_J dV=
           \rho_0 \int_{V_{\Box}} N_I \left( \bm{\xi} \right) N_J \left( \bm{\xi} \right) 
           J\left( \bm{\xi}; X_{Kn} \right) d\Box \\
           \left( I,J,K=1,\ldots,n_{nodes}, n=1,\ldots,n_{dim} \right)
    \end{gathered}
\end{equation}
Several lumping procedures are available in the finite element literature (see, e.g., \cite{dhondt2004finite}). In this study, we consider the standard row-sum lumping procedure, which produces a diagonal mass matrix. The subsequent developments are intended to be independent of the particular scaling technique   
\begin{equation}
    \label{eq: lumped_def}
    \begin{gathered}
    M_{II}^{diag}=\sum_{J=1}^{n_{nodes}} M_{IJ}^{cons} =  \sum_{J=1}^{n_{nodes}} M_{JI}^{cons} = 
           \int_{V_e} \rho_0 N_IdV=
           \rho_0 \int_{V_{\Box}} N_I \left( \bm{\xi} \right)  
           J\left( \bm{\xi}; X_{K,n} \right) d\Box \\
           M_{IJ}^{diag}{\Big|}_{I \neq J}=0, \qquad \left( I,K=1,\ldots,n_{nodes}, n=1,\ldots,n_{dim} \right)
    \end{gathered}
\end{equation}
The mixed mass matrix is defined as a convex combination of the lumped and consistent mass matrices, controlled by the parameter \( \alpha \). The limiting cases \( \alpha=0 \) and \( \alpha=1 \) recover the lumped and consistent formulations, respectively. Intermediate values balance computational efficiency against the accuracy of the inertia representation. 
\begin{equation}
    \label{eq: mixed_def}
    \begin{gathered}
    M_{IJ}^{mixed}= \left( 1-\alpha \right) M_{IJ}^{diag}+ \alpha M_{IJ}^{cons}, \qquad
           \left( I,J=1,\ldots,n_{nodes} \right)
    \end{gathered}
\end{equation}
While a purely lumped mass matrix $\alpha = 0$ ensures high computational efficiency by allowing trivial inversion, it introduces numerical dispersion error. Conversely, the consistent mass matrix $\alpha = 1$ provides superior frequency accuracy and captures inter-nodal inertia coupling, but requires expensive linear-system solves. By treating $\alpha$ as a weighting or blending parameter (often within $[0, 1]$), analysts can fine-tune the eigenvalue spectrum \cite{Felippa2014}. \\
\noindent \indent Mass matrices are real, Symmetric, and Positive Definite (SPD matrices). Linear combinations (\ref{eq: mixed_def}) of mass matrices and mass matrix inverses preserve these properties (e.g., \cite{HornJohnson1985}). \\
\noindent \indent We next consider the numerical evaluation of element mass matrices using standard numerical integration. In this approach, the integrand is evaluated at a prescribed set of quadrature points, and the resulting values are multiplied by the corresponding integration weights. The number and distribution of quadrature points determine the integration scheme's polynomial exactness and computational cost. Applying this procedure to the consistent mass matrix yields the approximation in Eq.~(\ref{eq: num_int_CMM}). We can also rewrite the expression as a linear combination of geometry-independent coefficient matrices, each multiplied by the Jacobian determinant evaluated at a quadrature point, as shown in Eq.~(\ref{eq: ST_CMM}). 
\begin{equation}
    \label{eq: num_int_CMM}
    \begin{gathered}
    M_{IJ}^{cons} \approx \rho_{0} \sum_{q=1}^{n_{ip}} w_q N_{Iq} N_{Jq} J_q \left( X_{Kn} \right), \qquad
    \left( I,J,K=1,\ldots,n_{nodes}, n=1,\ldots,n_{dim} \right)
    \end{gathered}
\end{equation}
where \( w_q \) are weights associated with each point, the number of integration points is denoted by \( n_{ip} \), while ~\( N_{Iq}={N_I}{\Big|}_{point \; q} \) and ~\( J_q={J}{\Big|}_{point \; q} \), \( \left( q=1,\ldots,n_{ip} \right) \). 
\begin{equation}
    \label{eq: ST_CMM}
    \begin{gathered}
    M_{IJ}^{cons \; ST} = \rho_0 \left( M_{IJ1}^{c \; ST} J_1 + M_{IJ2}^{c \; ST} J_2 + ... + M_{IJ{n_{ip}}}^{c \; ST} J_{n_{ip}} \right) =
        \rho_0 \sum_{q=1}^{n_{ip}} M_{IJq}^{c \; ST} J_q \left( X_{Kn} \right) \\
        M_{IJq}^{c \; ST}=w_q N_{Iq} N_{Jq} , \qquad
        \left( I,J,K=1,\ldots,n_{nodes}, n=1,\ldots,n_{dim} \right)
    \end{gathered}
\end{equation}
\noindent \indent This representation makes the mass matrix's dependence on element geometry explicit. For polynomial integrands, a sufficiently accurate quadrature rule can integrate the expression exactly. In practice, however, the quadrature order is typically selected to provide adequate accuracy at an acceptable computational cost. \\
\noindent \indent This observation motivates the following question: can the coefficient matrices be constructed so that a substantially shorter series provides an accurate approximation of the consistent mass matrix? The SA approach \cite{HG2_2016_reduce_int_point,hanukah2018improving} and the subsequent JD method \cite{hanukah_2025_element,hanukah_2025_JD10tet} address this question.\\
\noindent \indent Jacobian determinant representation (\ref{eq: J_SA}) together with consistent mass matrix formulation (\ref{eq: cmm_def}) results in an SA approximation 
\begin{equation}
    \label{eq: SA_CMM}
    \begin{gathered}
    M_{IJ}^{cons \; SA} = \rho_0 \left( M_{IJ1}^{c \; SA} J_1 + M_{IJ2}^{c \; SA} J_2 + ... + M_{IJn}^{c \; SA} J_{n} \right) =
        \rho_0 \sum_{q=1}^{n} M_{IJq}^{c \; SA} J_q \left( X_{Kn} \right) \\
        M_{IJq}^{c \; SA}= \int_{V_{\Box}} N_I N_J \hat{N}_q d\Box , \qquad
        \left( I,J,K=1,\ldots,n_{nodes}, n=1,\ldots,n_{dim} \right)
    \end{gathered}
\end{equation}
where coefficient matrices \( M_{IJq}^{c \; SA} \), \( \left( I, J=1,\ldots,n_{nodes},q=1,\ldots,n \right) \) are precomputed real symmetric matrices; since the integrand in the above is polynomial (multiplication if polynomials), numerical values of all the entries are evaluated using any Computer Algebra System (CAS) prior to in-code implementation. Practically, for a \( n_{dim}=3 \) element, four sampling points, i.e., a four-term series (\ref{eq: SA_CMM}), are sufficiently accurate. \\
\noindent \indent The JD representation of the Jacobian determinant (\ref{eq: J_JD}) together with the mass matrix definition (\ref{eq: cmm_def}) yields the JD approximation (representation)
\begin{equation}
    \label{eq: JD_CMM}
    \begin{gathered}
    M_{IJ}^{cons \; JD} = \rho_0 \left( M_{IJ0}^{c \; JD} \tilde{J}_0 + M_{IJ1}^{c \; JD} \tilde{J}_1 + ... + M_{IJ(n-1)}^{c \; JD} \tilde{J}_{n-1} \right)=
        \rho_0 \sum_{q=0}^{n-1} M_{IJq}^{c \; JD} \tilde{J}_q \left( X_{Kn} \right) \\
        M_{IJq}^{c \; JD}= \int_{V_{\Box}} N_I N_J \tilde{N}_q d\Box , \qquad
        \left( I,J,K=1,\ldots,n_{nodes}, n=1,\ldots,n_{dim} \right)
    \end{gathered}
\end{equation}
Here, coefficient matrices \( M_{IJq}^{c \; JD} \) are also precomputed weights, real and symmetric. Partial derivatives \( \tilde{J}_q \) are mesh-dependent. For the complete order polynomials \( \tilde{N}_q \) in (\ref{eq: J_JD}) and \( \hat{N}_q \) in (\ref{eq: J_SA}), the accuracies in (\ref{eq: SA_CMM}) and (\ref{eq: JD_CMM})  are practically identical and sufficient; however, for this accuracy level, the latter integrator (\ref{eq: JD_CMM}) is significantly more efficient. \\
\noindent \indent Let us adopt the SA and JD approaches to the lumped mass matrix (\ref{eq: lumped_def}), in order to expand it to the mixed (\ref{eq: mixed_def}) form 
\begin{equation}
    \label{eq: lumped_SA_JD}
    \begin{gathered}
        M_{IJ}^{diag \; SA} = \rho_0 \sum_{q=1}^{n} M_{IJq}^{d \; SA} J_q \left( X_{Kn} \right), \qquad 
        M_{IIq}^{d \; SA}= \int_{V_{\Box}} N_I \hat{N}_q d\Box \\
        M_{IJ}^{diag \; JD} = \rho_0 \sum_{q=0}^{n-1} M_{IJq}^{d \; JD} \tilde{J}_q \left( X_{Kn} \right), \qquad 
        M_{IIq}^{d \; JD}= \int_{V_{\Box}} N_I \tilde{N}_q d\Box \\
        \left( I,J,K=1,\ldots,n_{nodes}, n=1,\ldots,n_{dim} \right)
    \end{gathered}
\end{equation}
All the off-diagonal entries vanish. With the help of (\ref{eq: SA_CMM}),(\ref{eq: lumped_SA_JD}), the mixed mass matrix (\ref{eq: mixed_def}) is approximated (represented) following the SA approach
\begin{equation}
    \label{eq: mixed_SA}
    \begin{gathered}
        M_{IJ}^{mixed \; SA} = \rho_0 \left( M_{IJ1}^{m\;SA} J_1 + M_{IJ2}^{m\;SA} J_2 + ... + M_{IJn}^{m\;SA} J_n \right)=
        \rho_0 \sum_{q=1}^{n} M_{IJq}^{m \; SA} J_q \left( X_{Kn} \right) \\
        M_{IJq}^{m\;SA} = \left( 1-\alpha \right) M_{IJq}^{d \; SA}+\alpha M_{IJq}^{SA}, \qquad
        \left( I,J,K=1,\ldots,n_{nodes}, n=1,\ldots,n_{dim} \right)
    \end{gathered}
\end{equation}
With the help of (\ref{eq: JD_CMM}),(\ref{eq: lumped_SA_JD}), the mixed mass matrix (\ref{eq: mixed_def}) is approximated (represented) following the JD method
\begin{equation}
    \label{eq: mixed_JD}
    \begin{gathered}
        M_{IJ}^{mixed \; JD}=\rho_0 \left( M_{IJ0}^{m\;JD}\tilde{J}_0+M_{IJ1}^{m\;JD}\tilde{J}_1+...+M_{IJ(n-1)}^{m\;JD}\tilde{J}_{n-1} \right) =
        \rho_0\sum_{q=0}^{n-1} M_{IJq}^{m \; JD} \tilde{J}_q \left( X_{Kn} \right) \\
        M_{IJq}^{m\;JD} = \left( 1-\alpha \right) M_{IJq}^{d \; JD}+\alpha M_{IJq}^{JD}, \qquad
        \left( I,J,K=1,\ldots,n_{nodes}, n=1,\ldots,n_{dim} \right)
    \end{gathered}
\end{equation}
\noindent \indent In summary, the consistent, lumped, and mixed element mass matrices can all be expressed as linear combinations of geometry-independent coefficient matrices multiplied by geometry-dependent quantities. In standard numerical integration and the SA formulation, these quantities are Jacobian determinant evaluations at selected sampling points. In the JD formulation, they are the Jacobian determinant and its partial derivatives evaluated at the element centroid. These representations provide a common framework for constructing and comparing different mass-matrix integration schemes. For three-dimensional elements, the SA and JD approaches require relatively few terms to provide a sufficiently accurate approximation, and the JD formulation offers the computational advantages reported in \cite{hanukah_2025_element,hanukah_2025_JD10tet}. \\

\section{Direct Mass Inverse (DMI)}
\label{Sec. DMI}

The element-level mass matrix can be expressed as a linear combination of coefficient matrices of dimension ~\( n_{nodes} \times n_{nodes} \), weighted either by the Jacobian determinants evaluated at the sampling points, \( J_q \), as in (\ref{eq: mixed_SA}), or by the Jacobian determinants and their derivatives evaluated at the element centroid, \( \tilde{J}_q \), as in (\ref{eq: mixed_JD}). The exact inverse of the mixed element-level mass matrix is given by
\[  
W_{IJ}^{mixed} 
= 
\frac{ \operatorname{adj}( M_{KL}^{mixed} )_{IJ} }{ \det( M_{KL}^{mixed} ) }, 
\qquad 
(I,J,K,L=1,\ldots,n_{nodes})
\]
where \( \operatorname{adj}( \ldots )_{IJ} \) denotes the adjugate and \( \det(\ldots) \) denotes the determinant of the corresponding matrix. \\
\noindent \indent Provided that the entries of the mass matrix depend linearly on the Jacobian-related variables, the numerator, \( \operatorname{adj}( M_{KL}^{mixed} )_{IJ} \), is a matrix-valued polynomial of degree at most \( (n_{nodes}-1) \), whereas the denominator, \( \det(M_{KL}^{mixed}) \), is a polynomial of degree at most \( (n_{nodes}) \) in these variables. Accordingly, each component of the inverse mass matrix admits the rational representations
\begin{equation}
    \label{eq: W_exact}
    \begin{gathered}
W_{IJ}^{mixed} = 
\frac
{ \hat{P}_{IJ}^{mixed \; (n_{nodes}-1)} (J_1,\ldots,J_n) }
{ \hat{Q}^{mixed \; (n_{nodes})}   (J_1,\ldots,J_n) }
=
\frac
{ \tilde{P}_{IJ}^{mixed \; (n_{nodes}-1)} (\tilde{J}_0,\ldots,\tilde{J}_{n-1}) }
{ \tilde{Q}^{mixed \; (n_{nodes})}   (\tilde{J}_0,\ldots,\tilde{J}_{n-1}) },
\qquad (I,J=1,\ldots,n_{nodes}) 
    \end{gathered}
\end{equation}
where \( \hat{P}_{IJ}^{mixed} = \operatorname{adj}( M_{KL}^{mixed \; SA} )_{IJ} \), \( \hat{Q}^{mixed}=\det( M_{KL}^{mixed \; SA} ) \), and \( \tilde{P}_{IJ}^{mixed} = \operatorname{adj}( M_{KL}^{mixed \; JD} )_{IJ} \), \( \tilde{Q}^{mixed}=\det( M_{KL}^{mixed \; JD} ) \), denote the corresponding polynomial functions of the Jacobian-related variables. \\
\noindent \indent Although these expressions (\ref{eq: W_exact}) can, in principle, be derived and manipulated symbolically using a Computer Algebra System (CAS), their explicit evaluation and subsequent implementation are computationally impractical. In particular, the resulting expressions may become prohibitively large, making their direct incorporation into a numerical finite element implementation inefficient and difficult to maintain.\\
\noindent \indent Having established the exact rational structure of the inverse mass matrix, the central objective is therefore to construct an approximation that is computationally tractable, admits systematic refinement, and recovers the exact solution for the simplest admissible geometry: a regular element for which the Jacobian determinant is constant, such that \( J=J_q=\tilde{J}_0=J{\big|}_{\bm{\xi}^*}=const \) where \( \bm{\xi}^* \) denotes the element centroid. \\

\noindent \indent \textbf{Approximation strategy}. A natural starting point is the multivariate Taylor expansion of the inverse mass matrix about the regular-element configuration. Let \( \{ J\}=( J_1,\ldots,J_n ) \) and \( \{ \tilde{J} \} = ( \tilde{J}_0,\ldots,\tilde{J}_{(n-1)} ) \) denote the vector of the Jacobian-related variables, such that for the nominal operating point, about which all the approximations will be made, \( \{ J_0\}=( J_1=J_0,J_2=J_0,\ldots,J_n=J_0 ) \) and \( \{ \tilde{J}_0 \} = ( \tilde{J}_0,\tilde{J}_1=0,\tilde{J}_2=0,\ldots,\tilde{J}_{(n-1)}=0 ) \). Each component of the inverse mass matrix can then be regarded as a function \( \hat{W}_{IJ}^{mixed}( \{ J \} ) \) or \( \tilde{W}_{IJ}^{mixed}( \{ \tilde{J} \} ) \). If this function is sufficiently smooth in a neighborhood of the regular reference mesh (geometry), its multivariate Taylor expansion about \( \{ J_0 \} \) or \( \{\tilde{J}_0\} \) is
\[
W_{IJ}^{mixed}(\{ J \})=W_{IJ}^{mixed}(\{ J_0 \})+\sum_{\beta=1}^{n} \frac{\partial W_{IJ}^{mixed}}{\partial J_{\beta}}{\Big|}_{\{ J_0 \}}(J_{\beta}-J_0) + \ldots
\]
Thus
\[
W_{IJ}^{mixed \; SA} =
\hat{W}_{IJ0}^{m \; SA} \frac{1}{\rho_0 J_0} + 
\hat{W}_{IJ1}^{m \; SA} \frac{(J_1-J_0)}{\rho_0 J_0^2} + \ldots + 
\hat{W}_{IJn}^{m \; SA} \frac{(J_n-J_0)}{\rho_0 J_0^2} +
\hat{W}_{IJ(n+1)}^{m \; SA} \frac{(J_1-J_0) (J_2-J_0)}{\rho_0 J_0^3} + H.O.T.
\]
By reformulating, we write
\begin{equation}
    \label{eq: W_taylor_SA}
    \begin{gathered}
    W_{IJ}^{mixed \; SA} = \frac{1}{\rho_0} (
    W_{IJ0}^{m \; SA} \frac{1}{J_0} + 
    W_{IJ1}^{m \; SA} \frac{J_1}{J_0^2} + \ldots + 
    W_{IJn}^{m \; SA} \frac{J_n}{J_0^2} +
    W_{IJ(n+1)}^{m \; SA} \frac{J_1 J_2}{J_0^3}+\ldots )
    \end{gathered}
\end{equation}
\begin{equation}
    \label{eq: W_taylor_JD}
    \begin{gathered}
    W_{IJ}^{mixed \; JD} = \frac{1}{\rho_0} (
    W_{IJ0}^{m \; JD} \frac{1}{\tilde{J}_0} + 
    W_{IJ1}^{m \; JD} \frac{\tilde{J}_1}{\tilde{J}_0^2} + \ldots + 
    W_{IJ(n-1)}^{m \; JD} \frac{\tilde{J}_{n-1}}{\tilde{J}_0^2} +
    W_{IJn}^{m \; JD} \frac{\tilde{J}_1 \tilde{J}_2}{\tilde{J}_0^3}+\ldots )
    \end{gathered}
\end{equation}
where \( H.O.T. \) stands for Higher-Order Terms. The Taylor expansion provides a systematic local approximation and recovers the exact inverse mass matrix at the regular reference configuration  \( \{ J_0 \} \) or \( \{ \tilde{J}_0 \} \). However, its accuracy away from that configuration depends on the truncation order and on the behavior of the higher-order terms. In addition, the number of multivariate monomials increases rapidly with the expansion order and the number of independent geometric variables. Although this does not necessarily preclude a Taylor-based approximation, it significantly limits its efficiency and accuracy when a broad range of element geometries must be represented with a relatively compact expression. \\
\noindent \indent Therefore, based on our recent findings, we suggest practical guidelines for applying (\ref{eq: W_taylor_SA}) and (\ref{eq: W_taylor_JD}). Only linear elements can be treated accurately and efficiently with the above formulas, and only the linear part of the approximation should be considered. Regarding the efficiency, notice that the linear part of (\ref{eq: W_taylor_SA}) and (\ref{eq: W_taylor_JD}) is practically identical in its mathematical form and the required computational operations to (\ref{eq: mixed_SA}) and (\ref{eq: mixed_JD}), as it is almost the same formula with distinct coefficient matrices. Namely, at roughly the same number of operations required for element-level mass matrix evaluation, we can directly evaluate the inverse of the element-level mass matrix. We stress it once again: there is no need to evaluate the mass matrix or to invoke expensive inverse algorithms. Regarding accuracy and the must-have positive definiteness. We consider a set of highly distorted, yet valid \(J_q>0\), elements. For each, we ensure that our formula is PD. Then, we compute an exact inverse for each element and assess the proposed formula's accuracy using eigenvalues or the overall accuracy of the matrix components. \\
\noindent \indent \textbf{Padé approximation} offers an alternative, representing a function as the ratio of two polynomials rather than a single polynomial. In the multivariate setting, a Padé-type approximation of the inverse mass matrix can be written as
\[
W_{IJ}^{mixed} \approx R_{IJ}^{(L,M)}=\frac{P_{IJ}^{(L)}}{Q_{IJ}^{(M)}} 
\]
where \( P_{IJ}^{(L)} \) and \( Q_{IJ}^{(M)} \) are polynomials of prescribed degrees \( L \) and \( M \), respectively. Their coefficients can be determined by matching the multivariate Taylor coefficients of the exact function to a selected order, subject to a normalization such as \( Q_{IJ0}^{M}=1\), \( for \; all \ (I, J)  \), where  \(  Q_{IJ0}^{M} \) is the first term in the denominator series. This construction follows the general principle of rational approximation and its multivariate extensions \cite{BakerGravesMorris1996, Cuyt2000}.
\noindent \indent The main motivation for adopting a Padé-type approximation in this formulation is the rational structure of the exact inverse mass matrix. Specifically, the adjugate-over-determinant representation in (\ref{eq: W_exact}) shows that each component of the inverse is inherently rational in the Jacobian-related variables. A polynomial Taylor approximation does not explicitly retain this structure, whereas a rational approximation can represent both numerator and denominator effects in a single compact expression. Consequently, for a prescribed approximation complexity, a suitably constructed Padé approximant may provide improved accuracy over a wider range of geometric perturbations than a Taylor polynomial of comparable complexity. This advantage is particularly relevant when the determinant in the exact expression varies appreciably with the element geometry. \\
\noindent \indent Systematic refinement can be pursued by increasing the polynomial degrees or the order of the matched expansion coefficients. Nevertheless, increasing the approximation order alone does not guarantee convergence for arbitrary multivariate Padé constructions. The admissible geometric domain must also be chosen so that the approximate denominator remains nonzero, and the accuracy and stability of the resulting expressions must be assessed over the intended range of element distortions. These considerations form the basis for developing a compact, systematically refinable approximation of the inverse element-level mass matrix. \\
\noindent \indent From a practical perspective, we recommend starting with Padé approximation \( R_{IJ}^{(1,1)} \) and continuing to at most \( R_{IJ}^{(2,2)} \) for the second-order elements. The resulting matrix should be SPD, sufficiently accurate, and efficient. Following the SA approximation of the mixed mass matrix (\ref{eq: mixed_SA}), the first-order Padé approximation \( R_{IJ}^{(1,1)} \) takes the form
\begin{equation}
    \label{eq: W_pade_sa}
    \begin{gathered}
    W_{IJ}^{mixed} \approx \frac
    { \frac{1}{\rho_0} (
    \hat{P}_{IJ0} \frac{1}{J_0} + \hat{P}_{IJ1} \frac{J_1}{J_0^2} + \ldots + \hat{P}_{IJn} \frac{J_n}{J_0^2}
    )}
    {\hat{Q}_{IJ0} + \hat{Q}_{IJ1} \frac{J_1}{J_0} + \ldots + \hat{Q}_{IJn} \frac{J_n}{J_0}}, 
    \qquad \hat{Q}_{IJ0}=1
    \end{gathered}
\end{equation}
Following the JD approximation of the mixed mass matrix (\ref{eq: mixed_JD}), the first-order Padé approximation \( R_{IJ}^{(1,1)} \) takes the form
\begin{equation}
    \label{eq: W_pade_JD}
    \begin{gathered}
    W_{IJ}^{mixed} \approx \frac
    { \frac{1}{\rho_0} (
    \tilde{P}_{IJ0} \frac{1}{\tilde{J}_0} + \tilde{P}_{IJ1} \frac{\tilde{J}_1}{\tilde{J}_0^2} + \ldots + \tilde{P}_{IJ(n-1)} \frac{\tilde{J}_{n-1}}{\tilde{J}_0^2}
    )}
    {\tilde{Q}_{IJ0} + \tilde{Q}_{IJ1} \frac{\tilde{J}_1}{\tilde{J}_0} + \ldots + \tilde{Q}_{IJ(n-1)} \frac{\tilde{J}_{n-1}}{\tilde{J}_0}}, 
    \qquad \tilde{Q}_{IJ0}=1 
    \end{gathered}
\end{equation}
\noindent \indent The zeroth-order Padé approximation $R_{IJ}^{(0,0)}$, where the numerator and denominator consist of only one term, with the normalization $\tilde{Q}_{IJ0}=1$, recovers the regular element geometry exactly and becomes identical to the one-term series (\ref{eq: W_taylor_SA}), (\ref{eq: W_taylor_JD}), so that $W_{IJ}^{mixed}=\frac{1}{\rho_0 \tilde{J}_0} \tilde{P}_{IJ0} $. \\
\noindent \indent An important consideration is that the linear Padé approximations in (\ref{eq: W_pade_sa}) and (\ref{eq: W_pade_JD}) may involve different numbers of terms because they are based on different representations of the mixed mass matrix, namely, the SA and JD formulations. Although (\ref{eq: mixed_SA}) and (\ref{eq: mixed_JD}) may involve up to $n$ Jacobian-related variables, significantly shorter expansions may be sufficient in practice to achieve the desired balance between accuracy and computational efficiency.\\
\noindent \indent Therefore, the number of terms retained in the approximation should be selected based on numerical assessment rather than solely on the maximum number of available geometric variables. The principal selection criteria are preserving positive definiteness, computational efficiency, and sufficient accuracy of the resulting matrix.\\
\noindent \indent The computational efficiency of the proposed linear Padé approximations in (\ref{eq: W_pade_sa}) and (\ref{eq: W_pade_JD}) can be estimated by examining their mathematical structure. Both the numerator and denominator involve expressions structurally similar to those arising in the corresponding mass-matrix integration formulations (\ref{eq: mixed_SA}) and (\ref{eq: mixed_JD}). Consequently, directly evaluating the approximate inverse requires operations of the same general type as those used for mass-matrix integration, plus additional operations to evaluate the denominator and form the final ratio.\\
\noindent \indent In this respect, the proposed approach avoids explicit mass-matrix formation and subsequent numerical inversion. Furthermore, the SA and JD formulations offer significant computational advantages over conventional ST mass integration.\\

\section{Regular geometry model}
\label{Sec. RGM}

\noindent \indent Here, we consider the widely used eight-node isoparametric hexahedral element; see, e.g., \cite{dhondt2004finite}. An element with regular geometry, for which the Jacobian determinant is constant, $J=\mathrm{const}$, is a parallelepiped. Specifically, when the element nodes are configured such that each pair of opposing faces is parallel, the Jacobian determinant becomes independent of the local coordinates $\bm{\xi}=(\xi,\eta,\zeta)$,  i.e., $J=J(X_{Kn})=\mathrm{const}$. In this case, both mass matrix integrators SA (\ref{eq: mixed_SA}) and JD (\ref{eq: mixed_JD}) yield the exact result
\[
M_{IJ}^{mixed}=M_{IJ0}^{m \; JD} \rho_0 \tilde{J}_0
\]
This mass matrix is a positive linear combination of two symmetric positive definite (SPD) matrices. Consequently, for every $\alpha>0$, provided that the Jacobian determinant at the element centroid is positive, $J>0$, the resulting mass matrix is SPD. Its exact inverse is given by
\[
W_{IJ}^{mixed}=(M_{IJ0}^{m \; JD})^{-1} \frac{1}{\rho_0 \tilde{J}_0}
\]
which is also SPD. Importantly, all our proposed direct inverse formulae, (\ref{eq: W_taylor_SA}), (\ref{eq: W_taylor_JD}), (\ref{eq: W_pade_sa}), and (\ref{eq: W_pade_JD}), reproduce the regular element geometry exactly.
Given the following definitions,
\begin{equation}
    \label{eq: components}
    \begin{gathered}
    d = 416 \alpha^3 - 1524 \alpha^2 + 1836 \alpha - 729 \\
    m_1 = 52 \alpha^3 - 654 \alpha^2 + 1323 \alpha - 729, \qquad
    m_2 = 52 \alpha^3 - 156 \alpha^2 + 108 \alpha \\
    m_3 = 52 \alpha^3 - 108 \alpha^2 + 54 \alpha, \qquad
    m_4 = 52 \alpha^3 - 78 \alpha^2 + 27 \alpha
    \end{gathered}
\end{equation}
the inverse of an element-level mixed mass matrix takes the following form:
\begin{equation}
    \label{eq: z_o_matrix}
    \begin{gathered}
    \frac{1}{\rho_0 \; \tilde{J}_0 \; d} 
\left[\begin{array}{cccccccc}
m_1  & m_2  & m_3  & m_2  & m_2  & m_3  & m_4  & m_3  
\\
 m_2  & m_1  & m_2  & m_3  & m_3  & m_2  & m_3  & m_4  
\\
 m_3  & m_2  & m_1  & m_2  & m_4  & m_3  & m_2  & m_3  
\\
 m_2  & m_3  & m_2  & m_1  & m_3  & m_4  & m_3  & m_2  
\\
 m_2  & m_3  & m_4  & m_3  & m_1  & m_2  & m_3  & m_2  
\\
 m_3  & m_2  & m_3  & m_4  & m_2  & m_1  & m_2  & m_3  
\\
 m_4  & m_3  & m_2  & m_3  & m_3  & m_2  & m_1  & m_2  
\\
 m_3  & m_4  & m_3  & m_2  & m_2  & m_3  & m_2  & m_1  
\end{array}\right]
    \end{gathered}
\end{equation}
\noindent \indent The above expression provides an exact, non-diagonal inverse for elements with regular geometry, while serving as an approximation for the general case, $J\neq\mathrm{const}$.  Moreover, it is positive definite for all $\alpha>0$ and $J>0$. 
An important feature of this formulation is its robustness with respect to element distortion. Even when the Jacobian determinant vanishes at certain points within the element domain, while remaining positive at the centroid, the proposed expression remains well-defined and positive definite. Thus, it can potentially prevent the global algorithm from breaking down due to local geometric distortions. \\
\noindent \indent Furthermore, the parameter $\alpha$ provides an additional degree of freedom for calibrating the formulation and identifying its best practical value. The proposed approach is physically motivated, as it exactly reproduces the inverse mass matrix for regular element geometries, thereby preserving the consistency property associated with full integration. \\
\noindent \indent From a computational standpoint, the proposed expression (\ref{eq: z_o_matrix}) has a mathematical structure and computational cost comparable to those of the one-point ST integrator (\ref{eq: ST_CMM}) and the one-term series approximations in the SA and JD integrators (\ref{eq: mixed_SA}) and (\ref{eq: mixed_JD}). This is particularly attractive because it provides a direct, non-diagonal approximation of the inverse mass matrix at a computational cost comparable to that of the least expensive numerical integration schemes, without requiring a separate matrix inversion algorithm. \\
\noindent \indent Consequently, the proposed formula can be implemented selectively, on an element-by-element basis, in regions of coarse or highly distorted meshes, as well as for complete simulations. \\

\section{Conclusions}
\label{Sec. conclusions}

\noindent \indent This paper presents a systematic approach to directly evaluating the inverse of element-level mass matrices in solid isoparametric finite elements. The proposed method addresses mixed mass matrices defined as linear combinations of the lumped and consistent mass matrices, unifying diagonal, consistent, and intermediate mass formulations. \\
\noindent \indent The study yields two principal results. First, we develop a general, systematic formulation to directly evaluate the inverse of the element-level mixed mass matrix. Each entry of the general inverse is expressed as a polynomial in Jacobian-related variables divided by another polynomial in the same variables. The accuracy of the resulting approximation can be systematically improved by adding terms, and convergence to the exact inverse is established. The formulation applies to standard solid isoparametric elements in two and three dimensions, including lower- and higher-order elements. We propose alternative procedures for determining the coefficients, based on analytical derivation. Second, we derive a detailed closed-form formula for the inverse of the mixed mass matrix of the widely used eight-node hexahedral element. We examine the properties of this formula, with particular attention to its computational efficiency and numerical robustness. Notably, the formulation preserves positive definiteness when the Jacobian determinant at the element centroid is positive, including cases with coarse or highly distorted elements. In terms of computational cost, the proposed formula requires far fewer operations than conventional approaches that involve numerical integration followed by numerical matrix inversion. In particular, it costs less than standard one-point numerical integration of the mass matrix, even without accounting for the subsequent inversion. \\
\noindent \indent The proposed method is particularly relevant to explicit transient finite element analysis, where lumped mass matrices are widely used to avoid solving global systems of equations for nodal accelerations. By enabling direct evaluation of non-diagonal element-level mass-matrix inverses, the method provides a computationally efficient way to incorporate more representative inertial descriptions into explicit finite element formulations. This is particularly relevant to coarse- and distorted-mesh simulations, in which the limitations of diagonal mass approximations become more pronounced. The formulation can be applied selectively to coarse or distorted elements or across the entire mesh, balancing computational cost and solution accuracy. Its element-level structure also lends itself naturally to implementation within element-by-element computational procedures commonly employed in explicit finite element solvers. \\
\noindent \indent Future work will incorporate the proposed formulation into complete explicit transient finite element algorithms and assess its impact on the accuracy and computational cost of dynamic simulations through appropriate numerical benchmarks. Particular attention will be given to the influence of the mass-matrix weighting parameter, the selective application of non-diagonal mass formulations to coarse or distorted elements, and the method's performance across different element types and dynamic loading conditions. These investigations will establish the practical benefits of the proposed approach in engineering applications and further explore its potential for improving the accuracy and robustness of explicit solid dynamics simulations. \\

\noindent \textbf{Declaration of Generative AI and AI-Assisted Technologies} \\
During the preparation of this manuscript, the authors used Grammarly and ChatGPT (OpenAI) solely to improve the language, readability, and grammatical correctness of the text. The authors used these tools strictly for copy-editing, proofreading, and language refinement. The authors remain fully responsible for the technical content, engineering methodologies, and final conclusions presented in this paper. \\

\noindent \textbf{Acknowledgments} \\ 
Finally, I thank Ruth, Avital, and Orly for their endless support, understanding, and unwavering belief in my work, which gave me the peace of mind needed to complete this research. \\

\end{document}